\documentclass[10pt]{article}
\font\smallit=cmti10

\usepackage{amssymb,latexsym,amsmath,epsfig,amsthm} 

\makeatletter

\renewcommand\section{\@startsection {section}{1}{\z@}
{-30pt \@plus -1ex \@minus -.2ex}
{2.3ex \@plus.2ex}
{\normalfont\normalsize\bfseries\boldmath}}

\renewcommand\subsection{\@startsection{subsection}{2}{\z@}
{-3.25ex\@plus -1ex \@minus -.2ex}
{1.5ex \@plus .2ex}
{\normalfont\normalsize\bfseries\boldmath}}

\renewcommand{\@seccntformat}[1]{\csname the#1\endcsname. }

\makeatother

\newtheorem{theorem}{Theorem}

\theoremstyle{definition}
\newtheorem{definition}{Definition}

\newtheorem{remark}{Remark}
\newtheorem{example}{Example}

\begin{document}

\begin{center}
\uppercase{\bf A Method for Finding All Permutiples with a Fixed Set of Digits from a Single Known Example}
\vskip 20pt
{\bf Benjamin V. Holt}\\
{\smallit Department of Mathematics, Southwestern Oregon Community College, Oregon, USA}\\
{\tt benjamin.holt@socc.edu}\\
\vskip 10pt
\end{center}
\vskip 20pt

\centerline{\bf Abstract}

\noindent
A permutiple is a natural number that is a nontrivial multiple of a permutation of its digits in some base.
Special cases of permutiples include cyclic numbers (multiples of cyclic permutations of their digits) and palintiple numbers (multiples of their digit reversals). While cyclic numbers have a fairly straightforward description, palintiple numbers admit many varieties and cases. A previous paper attempts to get a better handle on the general case by constructing new examples of permutiples with the same set of digits, multiplier, and length as a known example. However, the results are not sufficient for finding all possible examples except when the multiplier divides the base. Using an approach based on the methods of this previous paper, we develop a new method which enables us to find all examples under any conditions. 

\pagestyle{myheadings}
\thispagestyle{empty}
\baselineskip=12.875pt
\vskip 30pt

\section{Introduction}

A {\it permutiple}, as the name suggests, is a number which is an integer multiple of some permutation of its digits in some natural number base, $b,$ greater than one \cite{holt_3}. Well-studied cases include cyclic numbers \cite{guttman,kalman}, that is, numbers which are multiples of cyclic permutations of their digits. A base-10 example of a cyclic number is $714285 = 5 \cdot 142857.$ 
A richer, but much less well-understood, case is palintiple numbers \cite{holt_1,holt_2}, also known as reverse multiples \cite{kendrick_1,sloane, young_1,young_2}, which are multiples of their digit reversals. The large variety of palintiple types can be organized using a graph-theoretical construction by Sloane \cite{sloane} called {\it Young graphs} which are a modification of the work of Young \cite{young_1,young_2}. The most widely known examples of palintiples, also in base 10, include $87912=4\cdot 21978$ and $98901 = 9 \cdot 10989.$ 

In \cite{holt_3}, the author establishes methods for finding new examples of general permutiples from old examples. For instance, using these methods, we are able to find new examples, such as $79128=4\cdot 19782,$ from the example with the same digits mentioned above. Although these methods shed some light on the problem, they are not able to account for all the desired examples under more general conditions. In particular, the results are not sufficient for finding the permutiple $78912=4\cdot 19728$ from the known example already mentioned. In this paper, we fill this gap by providing a simpler, yet more general, method for finding all permutiples with the same set of digits, multiplier, and length from a single known example.

\section{Basic Notation, Definitions, and Results}	

We shall use $(d_k,d_{k-1},\ldots,d_0)_b$ to denote the 
natural number $\sum_{j=0}^{k}d_j b^j$ where $0\leq d_j<b$ for all $0\leq j \leq k.$ The following is a definition of permutiple numbers.

\begin{definition}[\cite{holt_3}]
Let $n$ be a natural number and $\sigma$ be a permutation on the set $\{0,1,2,\ldots, k\}$.
We say that $(d_k, d_{k-1},\ldots, d_0)_b$  is an $(n,b,\sigma)$-\textit{permutiple} provided 
\[
(d_k,d_{k-1},\ldots,d_1, d_0)_b=n(d_{\sigma(k)},d_{\sigma(k-1)},\ldots, d_{\sigma(1)}, d_{\sigma(0)})_b.
\]
\end{definition}

For completeness, we recall a description of single-digit multiplication found in \cite{holt_3}. Letting $p_j$ denote the $j$th digit of the product, $c_j$ the $j$th carry, and $q_{j}$ the $j$th digit of the number being multiplied by $n$, we may perform single-digit, base-$b$ multiplication as follows:
\[
\begin{split}
c_0&=0,\\
p_j&=\lambda(n q_{j}+c_j),\\
c_{j+1}&=\left[ nq_{j}+c_j-\lambda(nq_{j}+c_j) \right]\div b,
\end{split}
\]
where $\lambda$ gives the least non-negative residue modulo $b$. Letting
$(p_k, p_{k-1}, \ldots, p_0)_b$ be the $(k+1)$-digit product, we have $c_{k+1}=0.$
Now, for a $(n,b,\sigma)$-permutiple, $(d_k, d_{k-1},\ldots, d_0)_b$, the above requires that $q_j=d_{\sigma(j)}$ and $d_{j}=p_j=\lambda(nd_{\sigma(j)}+c_j)$. What follows is the first result of \cite{holt_3}.

\begin{theorem}[\cite{holt_3}] \label{fund}
Let $(d_k, d_{k-1},\ldots, d_0)_b$ be an $(n,b,\sigma)$-permutiple and let $c_j$ be the
$j$th carry. Then
\[
b c_{j+1}-c_j=nd_{\sigma(j)}-d_{j}\]
for all $0\leq j \leq k$. 
\end{theorem}

Another important result from \cite{holt_3}, which will apply to our present purpose, is that the carries of any permutiple are less than the multiplier. We state this result formally.

\begin{theorem}[\cite{holt_3}] \label{carries}
Let $(d_k, d_{k-1},\ldots, d_0)_b$ be an $(n,b,\sigma)$-permutiple and let $c_j$ be the
$j$th
carry. Then $c_j\leq n-1$ for all $0 \leq j \leq k$.
\end{theorem}

Letting $\psi$ denote the $(k+1)$-cycle $(0,1,2,\ldots,k)$, we may write the above relations more conveniently in matrix form as
\[
(b P_\psi-I)\textbf{c} = (nP_\sigma-I) \textbf{d},
\]
where $I$ is the identity matrix, $P_\psi$ and $P_\sigma$ are permutation matrices, $\textbf{c}$ is a column vector containing the carries, and $\textbf{d}$ is a column vector containing the digits. We also note that indexing is from 0 to $k$ rather than from $1$ to $k.$ Finally, according to our description of single-digit multiplication, the first entry, $c_0$, of ${\bf c}$ is zero by definition.

The problem posed in \cite{holt_3} is the following: if $(d_k,d_{k-1},\ldots, d_0)_b$ is an $(n,b,\sigma)$-permutiple, find all permutations, $\pi,$ such that $(d_{\pi(k)},d_{\pi(k-1)},\ldots, d_{\pi(0)})_b$ is also a permutiple.

To sort through the types of new examples that arise, the notion of permutiple {\it conjugacy} was defined in \cite{holt_3}. Again, for completeness, we state this definition here.

\begin{definition}[\cite{holt_3}] \label{conj_class_def}
 Suppose $(d_k, d_{k-1},\ldots, d_0)_b$ is an $(n,b,\sigma)$-permutiple. Then an $(n,b, \tau_1)$-permutiple, $(d_{\pi_1(k)}, d_{\pi_1(k-1)},\ldots, d_{\pi_1(0)})_b,$ and
 an $(n,b, \tau_2)$-permutiple, $(d_{\pi_2(k)}, d_{\pi_2(k-1)},\ldots, d_{\pi_2(0)})_b$,
  are said to be \textit{conjugate} if $\pi_1 \tau_1 \pi_1^{-1}=\pi_2 \tau_2 \pi_2^{-1}$.
\end{definition}

Conjugacy defines an equivalence relation on the collection of permutiples having digits $d_k,$ $d_{k-1},$ $\ldots,$ $d_0.$ In \cite{holt_3}, these equivalence classes are refered to as {\it conjugacy classes}. The common permutation of a conjugacy class, $\beta=\pi_1 \tau_1 \pi_1^{-1}=\pi_2 \tau_2 \pi_2^{-1},$ is referred to as its {\it base permutation}. The methods of \cite{holt_3} are sufficient for finding all known examples within a conjugacy class, but fall short when trying to find all conjugacy classes outside of that which contains the known example.

\section{A Method for Finding All Examples}

We consider a generic $(n,b,\sigma)$-permutiple, $(d_k,d_{k-1},\ldots, d_0)_b,$ with carries $c_k,$ $c_{k-1},$ $\ldots,$ $c_0,$ and an $(n,b,\tau)$-permutiple with the same digits, $(d_{\pi(k)},d_{\pi(k-1)},\ldots, d_{\pi(0)})_b,$ but not necessarily the same carries,
$\hat{c}_k,$ $\hat{c}_{k-1},$ $\ldots,$ $\hat{c}_0.$
Then, in the notation established above, 
\begin{equation}
(nP_{\tau}-I)P_{\pi}{\bf d}=(bP_{\psi}-I){\bf \hat{c}}.
\label{fund_1}
\end{equation}
Reducing modulo $b,$ we have
\[
(nP_{\tau}-I)P_{\pi}{\bf d}\equiv -{\bf \hat{c}} \pmod{b}.
\]
Multiplying the above by $P_{\pi^{-1}}$ and rearranging, we obtain
\begin{equation}
{\bf d}+(b-n)P_{\pi \tau \pi^{-1}}{\bf d}\equiv P_{\pi^{-1}}{\bf \hat{c}} \pmod{b}.
\label{fund_2}
\end{equation}
Equation (\ref{fund_2}) in component form along with Theorem \ref{carries} give us our main result.

\begin{theorem} \label{main}
Suppose $(d_k,d_{k-1},\ldots, d_0)_b$ is an $(n,b,\sigma)$-permutiple. Then, in order for the number $(d_{\pi(k)},d_{\pi(k-1)},\ldots, d_{\pi(0)})_b$ to be an $(n,b,\tau)$-permutiple, it must be that 
\[
\lambda \left(d_j+(b-n)d_{\pi \tau \pi^{-1}(j)}\right)\leq n-1
\]
for all $0\leq j \leq k,$ where $\lambda$ is the least non-negative residue modulo $b.$
\end{theorem}
The above enables us to find all possible base permutations,
$\beta =\pi \tau \pi^{-1},$ for each conjugacy class by imposing necessary conditions on what $\pi \tau \pi^{-1}$ can be. A big advantage of the result is that it does not require any prior knowledge of what the carry sequence should be. In fact, once we narrow down the possible candidates for $\beta,$ we may then determine the values of the candidate set of carries by substituting in the known digits into Equation (\ref{fund_2}); the permuted carries are contained in the column vector ${\bf v}=P_{\pi^{-1}}{\bf \hat{c}}.$ 

With the base permutations in hand, we then rewrite Equation (\ref{fund_1}) as
\[
(nP_{\tau}-I)P_{\pi}{\bf d}=(bP_{\psi}-I)P_{\pi}{\bf v}
\]
since ${\bf \hat{c}}=P_{\pi}{\bf v}.$ Multiplying both sides by $P_{\pi^{-1}},$ we have
\[
(nP_{\pi \tau \pi^{-1}}-I){\bf d}=(bP_{\pi \psi \pi^{-1}}-I){\bf v},
\]
or,
\[
(nP_{\beta}-I){\bf d}=(bP_{\pi \psi \pi^{-1}}-I){\bf v}.
\]
Rearranging, we obtain
\begin{equation}
 bP_{\pi \psi \pi^{-1}}{\bf v}=(nP_{\beta}-I){\bf d}+{\bf v}.
 \label{fund_3}
\end{equation}
Now, since ${\bf d},$ ${\bf v},$ and $\beta=\pi \tau \pi^{-1}$ are known, the only unknown in Equation (\ref{fund_3}) is $\pi.$ This is to say that Equation (\ref{fund_3}) gives us a list of candidate  permutations, $\pi,$ for which $(d_{\pi(k)},d_{\pi(k-1)},\ldots, d_{\pi(0)})_b$ is an $(n,b,\tau)$-permutiple. Moreover, we note that Equation (\ref{fund_3}) is equivalent to Equation (\ref{fund_1}). So, Theorem 3 in \cite{holt_3} guarantees that every permutation, $\pi,$ satisfying Equation (\ref{fund_3}) yields a permutiple so long as $\hat{c}_0=0.$
From there, determining $\tau$ itself is a matter of either computing $\tau=\pi^{-1} \beta \pi$ or dividing $(d_{\pi(k)},d_{\pi(k-1)},\ldots, d_{\pi(0)})_b$ by $n.$

We now illustrate the above method by resolving a case for which the techniques of \cite{holt_3} were insufficient for recovering all permutiples from a known example.

\begin{example}
We shall find all 5-digit, $(4,10,\tau)$-permutiples with the same digits as the base-$10$ example $87912=4\cdot 21978.$ In more general notation, we state our known example as
 $(8,7,9,1,2)_{10}=4 \cdot (2,1,9,7,8)_{10}$ so that 
 \[
 {\bf d}
 =\left[                   
\begin{array}{c}
d_0\\d_1\\d_2\\d_3\\d_4\\
\end{array}
\right]
 =\left[                   
\begin{array}{c}
2\\1\\9\\7\\8\\
\end{array}
\right].
\]

Theorem \ref{main} tells us that for all $0 \leq j \leq 4$ we must have
\[
 \lambda \left(d_j+6d_{\pi \tau \pi^{-1}(j)}\right)\leq 3.
\]
That is,
\[
\begin{array}{c}
\lambda \left(2+6d_{\pi \tau \pi^{-1}(0)}\right)\leq 3,\\
\lambda \left(1+6d_{\pi \tau \pi^{-1}(1)}\right)\leq 3,\\
\lambda \left(9+6d_{\pi \tau \pi^{-1}(2)}\right)\leq 3,\\
\lambda \left(7+6d_{\pi \tau \pi^{-1}(3)}\right)\leq 3,\\
\lambda \left(8+6d_{\pi \tau \pi^{-1}(4)}\right)\leq 3.\\
\end{array}
\]
The above inequalities yield the possibilities
\[
\pi\tau\pi^{-1}=\left(
\begin{array}{ccccc}
0 & 1 & 2 & 3 & 4\\
4 & 0\mbox{ or }3 & 0,2,\mbox{ or }3 & 1\mbox{ or }2 & 0,2,\mbox{ or }3\\
\end{array}
\right).
\]
We note that $\pi\tau\pi^{-1}(3)=2$ would give us a relation that is not a permutation, so we are left with
\[
\pi\tau\pi^{-1}=\left(
\begin{array}{ccccc}
0 & 1 & 2 & 3 & 4\\
4 & 0\mbox{ or }3 & 0,2,\mbox{ or }3 & 1 & 0,2,\mbox{ or }3\\
\end{array}
\right).
\]
Consequently, the following are the possible base permutations:
\[
\begin{array}{c}
\beta_1=\left(
\begin{array}{ccccc}
0 & 1 & 2 & 3 & 4\\
4 & 3 & 2 & 1 & 0\\
\end{array}
\right),\,
\beta_2=\left(
\begin{array}{ccccc}
0 & 1 & 2 & 3 & 4\\
4 & 0 & 3 & 1 & 2\\
\end{array}
\right),\\\\
\beta_3=\left(
\begin{array}{ccccc}
0 & 1 & 2 & 3 & 4\\
4 & 0 & 2 & 1 & 3\\
\end{array}
\right),\,
\beta_4=\left(
\begin{array}{ccccc}
0 & 1 & 2 & 3 & 4\\
4 & 3 & 0 & 1 & 2\\
\end{array}
\right).
\end{array}
\]

The reader will notice that $\beta_1$ is the reversal permutation, $\rho,$ and is the base permutation of our known example. Also, $\beta_1=\rho$ is the digit permutation appearing in our known example. It is here that we underscore, as in \cite{holt_3}, that a base permutation need not be a digit permutation itself in conjugacy classes outside the one which contains the known example. 

From here, we substitute ${\bf d}$ and each possible $\beta_j=\pi\tau\pi^{-1}$ into Equation (\ref{fund_2}) to determine ${\bf v}=P_{\pi^{-1}}{\bf \hat{c}},$ which gives a possible set of carries. To these, we then apply Equation (\ref{fund_3}) to recover $\pi.$

Applying Equation (\ref{fund_2}) to 
$\beta_1=\left(
\begin{array}{ccccc}
0 & 1 & 2 & 3 & 4\\
4 & 3 & 2 & 1 & 0\\
\end{array}
\right)$  gives 
\[
{\bf v}\equiv P_{\pi^{-1}}{\bf \hat{c}}\equiv 
\left[\begin{array}{c}
2\\1\\9\\7\\8
\end{array}\right]
+
6P_{\beta_1}
\left[\begin{array}{c}
2\\1\\9\\7\\8
\end{array}\right]
\equiv
\left[\begin{array}{c}
0\\3\\3\\3\\0
\end{array}\right]
\pmod{10}.
\]
Since $0\leq \hat{c}_j\leq 3$ for each $0\leq j \leq 4,$ we conclude by Theorem \ref{carries} that
\[
{\bf v}=P_{\pi^{-1}}{\bf \hat{c}}=
\left[\begin{array}{c}
0\\3\\3\\3\\0
\end{array}\right],
\]
which is no surprise since these are the carries, $c_j,$ of our known example. 

Applying Equation (\ref{fund_3}),
\[
10 P_{\pi\psi\pi^{-1}}
\left[\begin{array}{c}
0\\3\\3\\3\\0
\end{array}\right]
=(4P_{\beta_1}-I)
\left[\begin{array}{c}
2\\1\\9\\7\\8
\end{array}\right]+\left[\begin{array}{c}
0\\3\\3\\3\\0
\end{array}\right]
=10\left[\begin{array}{c}
3\\3\\3\\0\\0
\end{array}\right].
\]
The possibilities are then expressed as
\[
\pi\psi\pi^{-1}=\left(
\begin{array}{ccccc}
0 & 1 & 2 & 3 & 4\\
1,2,\mbox{ or }3 & 1,2,\mbox{ or }3 & 1,2,\mbox{ or }3 & 0\mbox{ or }4 & 0\mbox{ or }4\\
\end{array}
\right).
\]
Since the above must be a 5-cycle, our possibilities are reduced to 
\[
\pi\psi\pi^{-1}=\left(
\begin{array}{ccccc}
0                  & 1                & 2               & 3 & 4\\
1\mbox{ or }2 & 2\mbox{ or }3  & 1 \mbox{ or } 3 & 4 & 0 \\
\end{array}
\right).
\]
That is, $\pi\psi\pi^{-1}=(\pi(0),\pi(1),\pi(2),\pi(3),\pi(4))$ can be either $\psi=(0,1,2,3,4)$ or $(0,2,1,3,4)=(1,2)\psi(1,2).$  Thus, $\pi=\varepsilon$, the identity permutation, and $\pi=(1,2)$ both solve Equation (3). 
Moreover, since the first carry, $\hat{c}_0=v_{\pi(0)},$ must always be zero, the application of Equation (\ref{fund_3}) above gives that $\pi(0)$ can be either $0$ or $4.$ The possibility $\pi(0)=4$ gives us two more permutations: $\pi=\psi^4$ and $\pi=(1,2)\psi^4.$

The entire conjugacy class for 
$\beta_1=\rho$ is listed below.
\begin{center}
\begin{footnotesize}
 \begin{tabular}{|c|c|c|c|}
\hline
$(d_{\pi(4)},d_{\pi(3)},d_{\pi(2)},d_{\pi(1)},d_{\pi(0)})_{10}$ & $\pi$ & $\tau$ & $(\hat{c}_4,\hat{c}_3,\hat{c}_2,\hat{c}_1,\hat{c}_0)$\\\hline
$(8,7,9,1,2)_{10}$ & $\varepsilon$ & $\rho$ & $(0,3,3,3,0)$\\\hline
$(8,7,1,9,2)_{10}$ & $(1,2)$ & $(1,2)\rho(1,2)$ & $(0,3,3,3,0)$\\\hline
$(7,9,1,2,8)_{10}$ & $\psi^4$ & $\psi^{-4}\rho\psi^4$ & $(3,3,3,0,0)$\\\hline
$(7,1,9,2,8)_{10}$ & $(1,2)\psi^4$ & $\psi^{-4}(1,2)\rho(1,2)\psi^4$ & $(3,3,3,0,0)$\\\hline
\end{tabular}
\end{footnotesize}
\end{center}

\begin{remark}
 The reader will notice that Equation (\ref{fund_3}) does the same work as Corollary 2 in \cite{holt_3}. In fact, the above analysis is identical in form to that found in Example 3 in \cite{holt_3}.
\end{remark}

We now find the conjugacy class for $\beta_2=\left(
\begin{array}{ccccc}
0 & 1 & 2 & 3 & 4\\
4 & 0 & 3 & 1 & 2\\
\end{array}
\right).$ Again, by Equation (\ref{fund_2}),
\[
{\bf v}\equiv P_{\pi^{-1}}{\bf \hat{c}}\equiv 
\left[\begin{array}{c}
2\\1\\9\\7\\8
\end{array}\right]
+
6P_{\beta_2}
\left[\begin{array}{c}
2\\1\\9\\7\\8
\end{array}\right]
\equiv
\left[\begin{array}{c}
0\\3\\1\\3\\2
\end{array}\right]
\pmod{10}.
\]
In similar fashion to the above case, we may argue that
\[
{\bf v}=P_{\pi^{-1}}{\bf \hat{c}}=
\left[\begin{array}{c}
0\\3\\1\\3\\2
\end{array}\right].
\]
Using Equation (\ref{fund_3}), we have
\[
10 P_{\pi\psi\pi^{-1}}
\left[\begin{array}{c}
0\\3\\1\\3\\2
\end{array}\right]
=(4P_{\beta_2}-I)
\left[\begin{array}{c}
2\\1\\9\\7\\8
\end{array}\right]+\left[\begin{array}{c}
0\\3\\1\\3\\2
\end{array}\right]
=10\left[\begin{array}{c}
3\\1\\2\\0\\3
\end{array}\right],
\]
which allows for
\[
\pi\psi\pi^{-1}=\left(
\begin{array}{ccccc}
0 & 1 & 2 & 3 & 4\\
1\mbox{ or }3 & 2 & 4 & 0 & 3\mbox{ or }1\\
\end{array}
\right).
\]
Again, since the above must be a 5-cycle, we are left with a single possibility,
\[
\pi\psi\pi^{-1}=\left(
\begin{array}{ccccc}
0 & 1 & 2 & 3 & 4\\
1 & 2 & 4 & 0 & 3\\
\end{array}
\right).
\]
It follows that $\pi\psi\pi^{-1}=(\pi(0),\pi(1),\pi(2),\pi(3),\pi(4))=(0,1,2,4,3).$ Since Equation (\ref{fund_3}) only allows for $\pi(0)=0$ in this case, the only possible permutation is $\pi=(3,4).$

The conjugacy class for $\beta_2,$ therefore, consists of a single element, namely, the example $(7,8,9,1,2)_{10}=4\cdot (1,9,7,2,8)_{10},$ a $(4,10,\tau)$-permutiple with carry vector
\[
{\bf \hat{c}}=P_{\pi}{\bf v}=
\left[\begin{array}{c}
0\\3\\1\\2\\3
\end{array}\right],
\]
where 
$\tau=
\left(
\begin{array}{ccccc}
0 & 1 & 2 & 3 & 4\\
3 & 0 & 4 & 2 & 1\\
\end{array}
\right)
=(0,3,2,4,1)=\pi^{-1}\beta_2\pi.
$
\begin{remark}
The above conjugacy class consists of the example that the results in \cite{holt_3} could not account for. 
\end{remark}

Considering $\beta_3=\left(
\begin{array}{ccccc}
0 & 1 & 2 & 3 & 4\\
4 & 0 & 2 & 1 & 3\\
\end{array}
\right),$ another use of Equation (\ref{fund_2}) yields
\[
{\bf v}\equiv P_{\pi^{-1}}{\bf \hat{c}}\equiv 
\left[\begin{array}{c}
2\\1\\9\\7\\8
\end{array}\right]
+
6P_{\beta_3}
\left[\begin{array}{c}
2\\1\\9\\7\\8
\end{array}\right]
\equiv
\left[\begin{array}{c}
0\\3\\3\\3\\0
\end{array}\right]
\pmod{10}.
\]
Then, using
\[
{\bf v}=
\left[\begin{array}{c}
0\\3\\3\\3\\0
\end{array}\right],
\]
we employ Equation (\ref{fund_3}) to obtain
\[
10 P_{\pi\psi\pi^{-1}}
\left[\begin{array}{c}
0\\3\\3\\3\\0
\end{array}\right]
=(4P_{\beta_3}-I)
\left[\begin{array}{c}
2\\1\\9\\7\\8
\end{array}\right]+\left[\begin{array}{c}
0\\3\\3\\3\\0
\end{array}\right]
=10\left[\begin{array}{c}
3\\1\\3\\0\\2
\end{array}\right].
\]
Since there is no permutation, $\pi,$ which makes the above statement true, we conclude that the conjugacy class corresponding to $\beta_3$ is empty.
By a similar calculation, $\beta_4$ also yields no new examples.
With the above, we have found all $(4,10,\tau)$-permutiples with the same digits as our known example. 
\end{example}
 
We further demonstrate the above techniques by presenting a base-13 example with a repeated digit.
 
\begin{example}
We shall consider a known base-13 example, 
\[
(9,1,12,12,3,11)_{13}=7\cdot(1,3,12,12,11,9)_{13}.
\]
In this case we do not have a unique permutation, $\sigma,$ due to a repeated digit. We shall choose the simpler of the two options, 
\[
\sigma=\left(\begin{array}{cccccc}0 & 1 & 2 & 3 & 4 & 5\\5 & 0 & 2 & 3 & 1 & 4\\\end{array}\right),
\]
knowing that $\pi=(2,3)$ will appear in the analysis which follows.

Assembling and labeling our digits in a column-vector format,
\[
{\bf d}
=\left[\begin{array}{c}d_{0}\\d_{1}\\d_{2}\\d_{3}\\d_{4}\\d_{5}\\\end{array}\right]
=\left[\begin{array}{c}11\\3\\12\\12\\1\\9\\\end{array}\right],
\]
we now invoke Theorem \ref{main}. That is, any base permutation, $\beta=\pi \tau \pi^{-1},$ must satisfy 
\[
\begin{array}{c}
 \lambda\left(11+6d_{\pi \tau \pi^{-1}(0)}\right)\leq 6,\\
 \lambda\left(3+6d_{\pi \tau \pi^{-1}(1)}\right)\leq 6,\\
 \lambda\left(12+6d_{\pi \tau \pi^{-1}(2)}\right)\leq 6,\\
 \lambda\left(12+6d_{\pi \tau \pi^{-1}(3)}\right)\leq 6,\\
 \lambda\left(1+6d_{\pi \tau \pi^{-1}(4)}\right)\leq 6,\\
 \lambda\left(9+6d_{\pi \tau \pi^{-1}(5)}\right)\leq 6.\\
 \end{array}
\]
The above gives the following possibilities:
\[
\begin{array}{rcl}\beta(0)=\pi \tau \pi^{-1}(0) &= & 1, 2, 3, 4, \mbox{ or } 5, \\\beta(1)=\pi \tau \pi^{-1}(1) &= & 0 \mbox{ or }  5, \\\beta(2)=\pi \tau \pi^{-1}(2) &= & 0, 1, 2, 3, 4, \mbox{or }  5, \\\beta(3)=\pi \tau \pi^{-1}(3) &= & 0, 1, 2, 3, 4, \mbox{ or }  5, \\\beta(4)=\pi \tau \pi^{-1}(4) &= & 0, 1, \mbox{ or } 5, \\\beta(5)=\pi \tau \pi^{-1}(5) &= & 1, 2, 3, \mbox{ or } 4. \\\end{array}
\]
We note here that due to the volume of possibilities presented by this example, we resorted to writing computer code to determine all possible base permutation candidates. There were 78 candidates in total, many of which we were able to immediately rule out as none of the entries of the vector ${\bf v}=P_{\pi^{-1}}{\bf \hat{c}}$ were 0. Moreover, as seen in the previous example, other candidates did not allow for Equation (\ref{fund_3}) to be satisfied as the vectors $bP_{\pi \psi \pi^{-1}}{\bf v}$ and $(nP_{\beta}-I){\bf d}+{\bf v}$ did not contain the same collection of entries. This enabled us to rule out more base permutation candidates. The following are the six viable base permutations for the original example:
\[
\begin{array}{c}
\beta_{1}=\left(\begin{array}{cccccc}0 & 1 & 2 & 3 & 4 & 5\\5 & 0 & 2 & 3 & 1 & 4\\\end{array}\right),\,
\beta_{2}=\left(\begin{array}{cccccc}0 & 1 & 2 & 3 & 4 & 5\\5 & 0 & 3 & 2 & 1 & 4\\\end{array}\right),\\\\ 
\beta_{3}=\left(\begin{array}{cccccc}0 & 1 & 2 & 3 & 4 & 5\\3 & 5 & 2 & 0 & 1 & 4\\\end{array}\right),\, 
\beta_{4}=\left(\begin{array}{cccccc}0 & 1 & 2 & 3 & 4 & 5\\3 & 5 & 0 & 2 & 1 & 4\\\end{array}\right),\\\\
\beta_{5}=\left(\begin{array}{cccccc}0 & 1 & 2 & 3 & 4 & 5\\2 & 5 & 3 & 0 & 1 & 4\\\end{array}\right),\,
\beta_{6}=\left(\begin{array}{cccccc}0 & 1 & 2 & 3 & 4 & 5\\2 & 5 & 0 & 3 & 1 & 4\\\end{array}\right).
\end{array}
\]
The reader will notice that $\beta_1$ and $\beta_2$ are the two choices of permutation we had for the original example. We will now consider both simultaneously since applying Equations (\ref{fund_2}) and (\ref{fund_3}) results in the same set of equations. That is, for both $\beta=\beta_1$ and $\beta=\beta_2,$ Equation (\ref{fund_2}) gives
\[
 {\bf v} \equiv P_{\pi^{-1}}{\bf \hat{c}}\equiv\left[\begin{array}{c}11\\3\\12\\12\\1\\9\\\end{array}\right]+6P_{\beta}\left[\begin{array}{c}11\\3\\12\\12\\1\\9\\\end{array}\right]\equiv\left[\begin{array}{c}0\\4\\6\\6\\6\\2\\\end{array}\right] \pmod{13}
\]
so that by Theorem \ref{carries},
\[
 {\bf v}=\left[\begin{array}{c}0\\4\\6\\6\\6\\2\\\end{array}\right].
\]
An application of Equation (\ref{fund_3}) then gives
\[
13P_{\pi \psi \pi^{-1}}\left[\begin{array}{c}0\\4\\6\\6\\6\\2\\\end{array}\right]=(7P_{\beta}-I)\left[\begin{array}{c}11\\3\\12\\12\\1\\9\\\end{array}\right]+\left[\begin{array}{c}0\\4\\6\\6\\6\\2\\\end{array}\right]=13\left[\begin{array}{c}4\\6\\6\\6\\2\\0\\\end{array}\right].
\]
\begin{remark}
 If we consider $d_2=12$ and $d_3=12$ to be distinct digits, then $\beta_1$ and $\beta_2$ really do represent two distinct conjugacy classes.
\end{remark}

We now consider possibilities for $\pi \psi \pi^{-1}.$ From the above use of Equation (\ref{fund_3}), we have that
\[
\pi \psi \pi^{-1}=
\left(
\begin{array}{cccccc}0 & 1 & 2 & 3 & 4 & 5\\1 & 2, 3,\mbox{ or } 4 & 2, 3,\mbox{ or } 4 & 2, 3,\mbox{ or } 4 & 5 & 0\\\end{array}
\right).
\]
Since this must be a 6-cycle, $(\pi(0),\pi(1),\pi(2),\pi(3),\pi(4),\pi(5))$, the possibilities are pared down to $(0,1,2,3,4,5)$ and $(0,1,3,2,4,5).$
Also, since $\hat{c}_{0}=v_{\pi(0)}=0,$ the above application of Equation (\ref{fund_3}) requires that $\pi(0)=0.$
It follows that $\pi$ must be the identity permutation or the transposition $(2,3).$

Thus, our original example, $(d_5,d_4,d_3,d_2,d_1,d_0)_{13}=(9,1,12,12,3,11)_{13},$ is the only element in the conjugacy class corresponding to $\beta_1.$ If we consider $d_2=12$ and $d_3=12$ to be distinct from one another, then 
\[
(d_{5},d_{4},d_{2},d_{3},d_{1},d_{0})_{13}=(9,1,12,12,3,11)_{13} 
\]
is the solitary element of the conjugacy class corresponding to $\beta_2.$

Moving on to other base permutations, applying Equation (\ref{fund_2}) to $\beta_3$ gives 
\[
 {\bf v}=\left[\begin{array}{c}5\\5\\6\\0\\6\\2\\\end{array}\right]
\]
so that Equation (\ref{fund_3}) becomes 
\[
13P_{\pi \psi \pi^{-1}}\left[\begin{array}{c}5\\5\\6\\0\\6\\2\\\end{array}\right]=(7P_{\beta_{3}}-I)\left[\begin{array}{c}11\\3\\12\\12\\1\\9\\\end{array}\right]+\left[\begin{array}{c}5\\5\\6\\0\\6\\2\\\end{array}\right]=13\left[\begin{array}{c}6\\5\\6\\5\\2\\0\\\end{array}\right].
\]
Then 
\[
\pi \psi \pi^{-1}=
\left(
\begin{array}{cccccc}0 & 1 & 2 & 3 & 4 & 5\\2 \mbox{ or } 4 & 0\mbox{ or } 1 & 2 \mbox{ or } 4 & 0 \mbox{ or } 1 & 5 & 3\\\end{array}
\right),
\]
of which the only 6-cycle is $(0,2,4,5,3,1).$ Again, imposing the restriction that $\hat{c}_{0}=v_{\pi(0)}=0,$ the above application of Equation (\ref{fund_3}) requires that $\pi(0)=3.$ Hence,
$
\pi=
\left(
\begin{array}{cccccc}0 & 1 & 2 & 3 & 4 & 5\\ 3 & 1 & 0 & 2 & 4 & 5\\\end{array}
\right).
$
From the above, we obtain a new example,
\[
(9,1,12,11,3,12)_{13}=7\cdot (1,3,12,12,9,11)_{13}.
\]
We now apply Equation (\ref{fund_2}) to $\beta_4$ to obtain 
\[
{\bf v}=\left[\begin{array}{c}5\\5\\0\\6\\6\\2\\\end{array}\right].
\]
Equation (\ref{fund_3}) then becomes
\[
13P_{\pi \psi \pi^{-1}}\left[\begin{array}{c}5\\5\\0\\6\\6\\2\\\end{array}\right]=(7P_{\beta_{4}}-I)\left[\begin{array}{c}11\\3\\12\\12\\1\\9\\\end{array}\right]+\left[\begin{array}{c}5\\5\\0\\6\\6\\2\\\end{array}\right]=13\left[\begin{array}{c}6\\5\\5\\6\\2\\0\\\end{array}\right].
\]
The reader will note that in this case, the permuted carry vector, 
${\bf v},$ has the same collection of entries as that in the case of $\beta_3.$ In fact, the above equation yields a transposition, $\pi=(0,2),$ which produces the same numerical example as $\beta_3,$ 
\[
(d_{5},d_{4},d_{3},d_{0},d_{1},d_{2})_{13}=(9,1,12,11,3,12)_{13}.
\]
Again, this is a result of the repeated digit. We again emphasize that if we keep track of which repeated digit is which, treating them as distinct objects, then the above example would be considered new.
Finally, we note that applying Equations (\ref{fund_2}) and (\ref{fund_3}) to $\beta_5$ and $\beta_6$ yields the same set of equations as $\beta_3$ and $\beta_4,$ respectively. 

Not distinguishing between repeated digits, we have shown that the single new example above is the only other one with the same digits as the original.

\end{example}

\section{Summary of Method and Concluding Remarks}

To summarize the above method, Theorem \ref{main} provides a list of base permutation candidates. Trying all of these in Equation (\ref{fund_2}) gives us possibilities for what the carries can be by computing the permuted carry vector ${\bf v}=P_{\pi^{-1}}{\bf \hat{c}}.$ Inserting this information into Equation (\ref{fund_3}) then allows us to recover permutations, $\pi,$ which yield new permutiples.  

While the above method addresses the main question posed in \cite{holt_3}, we note that there are still plenty of questions that remain from the above considerations. Of particular interest are patterns or restrictions which may exist concerning permutation type and orders of base permutations, as well as the sizes of their corresponding conjugacy classes.

Other tractable lines of inquiry with the goal of finding new permutiples from old include understanding when ``derived'' permutiples are possible, that is, $(n,b,\sigma)$-permutiples whose truncated carry vector, $(c_k,c_{k-1},\ldots,c_1),$ is also a base-$n$ permutiple. An example of this phenomenon, mentioned in \cite{holt_3}, is the cyclic $(6,12,\psi^3)$-permutiple, $(10,3,5,1,8,6)_{12}=6\cdot(1,8,6,10,3,5)_{12}$, whose nonzero carries are the digits of the $(2,6)$-palintiple, $(4,3,5,1,2)_6=2 \cdot (2,1,5,3,4)_6.$ In other words, we ask: when is it possible to construct, or ``derive,'' a  new permutiple, say $(10,3,5,1,8,6)_{12},$ from a known permutiple, such as $(4,3,5,1,2)_6,$ by treating it as a carry vector, $(4,3,5,1,2,0)?$ We note that the less general case of derived palintiples is taken up in \cite{holt_2}. 

For other questions and results regarding the general permutiple problem, the reader is directed to \cite{holt_3}.

\vskip 20pt\noindent {\bf Acknowledgement.} The author gratefully acknowledges the contribution of the anonymous referee. We extend our sincere thanks for their time and helpful suggestions which substantially improved this work.

\end{document}